\documentclass{amsart}
\usepackage{amssymb}

\usepackage{a4}
\begin{document}
\newtheorem{theorem}{Theorem}[section]
\newtheorem{lemma}[theorem]{Lemma}
\newtheorem{remark}[theorem]{Remark}
\newtheorem{definition}[theorem]{Definition}
\newtheorem{corollary}[theorem]{Corollary}
\newtheorem{example}[theorem]{Example}
\def\BB{\mathcal{B}}
\def\id{\operatorname{Id}}
\def\ffrac#1#2{{\textstyle\frac{#1}{#2}}}
\def\Span{\operatorname{Span}}
\def\Range{\operatorname{Range}}
\def\Rank{\operatorname{Rank}}
\def\Gr{\operatorname{Gr}}
\def\OGr{\widetilde{\operatorname{Gr}}}
\makeatletter
 \renewcommand{\theequation}{%
  \thesection.\alph{equation}}
 \@addtoreset{equation}{section}
 \makeatother
\title[Manifolds which are Ivanov-Petrova or $k$-Stanilov]
{Manifolds which are Ivanov-Petrova or $k$-Stanilov}
\author{P. Gilkey, S. Nik\v cevi\'c, and V. Videv}
\begin{address}{PG: Mathematics Department, University of Oregon,
Eugene Or 97403 USA.\newline Email: {\it gilkey@darkwing.uoregon.edu}}
\end{address}
\begin{address}{SN: Mathematical Institute, Sanu,
Knez Mihailova 35, p.p. 367,
11001 Belgrade,
Yugoslavia.
Email: {\it stanan@mi.sanu.ac.yu}}\end{address}
%\begin{address}{GS: Faculty of Mathematics and Informatics, University
%          of Sofia,  5 James Bourchier Blvd, 1164 Sofia Bulgaria.
%          email: tba}\end{address}
\begin{address}{VV: Mathematics Department, Thracian University, University
    Campus, 6000 Stara Zagora, Bulgaria.
     email: videv@uni-sz.bg}\end{address}
\keywords{skew-symmetric curvature operator, higher order curvature operator, Jordan Stanilov manifolds,
Jordan Ivanov-Petrova manifolds, curvature homogeneity.
\newline \phantom{.....}2000 {\it Mathematics Subject Classification.} 53B20.}

\begin{abstract} We present some examples of curvature homogeneous pseudo-Riemannian manifolds which are $k$-spacelike Jordan
Stanilov.
\end{abstract}
\maketitle

\section{Introduction}\label{Sect-1}

In considering the spectral geometry of the Riemann curvature tensor, one studies when a
certain natural operator associated to the curvature has constant Jordan normal form on the
natural domain of definition. In this brief note, we consider two such operators -- the
skew-symmetric curvature operator $\mathcal{R}$ and a higher order generalization $\Theta$. We
begin by recalling some basic definitions.

\subsection{The algebraic context} Let
$\mathcal{V}:=(V,g,R)$ be a model space where
$V$ is a finite dimensional real vector space which is equipped with a non-degenerate inner product $g$ of
signature $(p,q)$, and where $R$ is an algebraic curvature tensor on $V$, i.e.
$R\in\otimes^4V^*$ satisfies the usual curvature symmetries:
\begin{eqnarray}
&&R(x,y,z,w)=R(z,w,x,y)=-R(y,x,z,w),\label{eqn-1.xx}\\
&&R(x,y,z,w)+R(y,z,x,w)+R(z,x,y,w)=0\label{eqn-1.yy}\,.
\end{eqnarray}
Let $R(x,y)$ be the associated curvature operator; it is characterized by the identity:
$$R(x,y,z,w)=g(R(x,y)z,w)\,.$$

Let $\OGr_{k,\pm}(V,g)$ (resp. $\Gr_{k,\pm}(V,g)$) be the Grassmannians of oriented (resp.
unoriented) spacelike ($+$) and timelike ($-$) $k$-planes in $V$. Let $\{e_1,e_2\}$ be an oriented orthonormal basis for $\pi\in\OGr_{k,\pm}(V,g)$. The
skew-symmetric curvature operator
$$\mathcal{R}(\pi):=R(e_1,e_2)$$
was introduced by Stanilov in 1990 -- see the discussion in
Ivanova and Stanilov \cite{refRIGS};
it is independent of the particular oriented orthonormal basis chosen for $\pi$. One says $\mathcal{V}$ has {\it constant
spacelike} (resp. {\it timelike}) {\it rank} $r$ if $\Rank\{\mathcal{R}(\pi)\}=r$ for any oriented spacelike (resp. timelike) $2$
plane $\pi$ of $V$.  
One says $\mathcal{V}$ is {\it spacelike} (resp. {\it timelike}) {\it Jordan Ivanov-Petrova}
if the Jordan normal form of
$\mathcal{R}$ is constant on $\OGr_{2,\pm}(V,g)$.
Clearly if $\mathcal{V}$ is spacelike (resp. timelike) Jordan Ivanov-Petrova, then $\mathcal{V}$
has constant spacelike (resp. timelike) rank. 

There is a higher order analogue which was introduced Stanilov \cite {Stx,Sty}. If
$\{e_1,...,e_k\}$ is an orthonormal basis for $\pi\in\Gr_{k,\pm}(M,g)$, then the {\it higher order
curvature operator} is defined by setting:
$$\Theta(\pi):=\sum_{i,j}R(e_i,e_j)R(e_i,e_j)\,.$$
This self-adjoint operator is, similarly, independent of the particular orthonormal basis chosen for $\pi$.
One says $\mathcal{V}$ is {\it $k-$spacelike} (resp. {\it $k$-timelike}) {\it Jordan Stanilov}
if the Jordan normal form of $\Theta$ is constant on $\Gr_{k,\pm}(V,g)$; see also
\cite{TzVi99} for further details. Up to a suitable normalizing factor,
$$\Theta(\pi)=\textstyle\int_{\xi\in\tilde Gr_2(\pi)}\mathcal{R}(\xi)^2d\xi$$
so the higher order curvature operator can be regarded as an average of the
square skew-symmetric curvature operator. It is necessary to square $\mathcal{R}$ to obtain a non-zero average since
$\mathcal{R}(\cdot)$ changes sign if the orientation of $\pi$ is reversed.

Note that one could in fact define $\mathcal{R}$ (resp. $\Theta$) for
any non-degenerate oriented $2$ plane (resp. non-degenerate unoriented $k$ plane); we shall
restrict ourselves to the spacelike and the timelike planes in the interests of simplicity.

\subsection{The geometric context} Let $(M,g)$ be a pseudo-Riemannian manifold of signature
$(p,q)$ and dimension $m=p+q$. Let $R$ be the Riemann curvature of the Levi-Civita
connection. We say that $(M,g)$ is {\it spacelike} (resp. {\it timelike}) Jordan
Ivanov-Petrova if $\mathcal{T}_P:=(T_PM,g_P,R_P)$ is spacelike (resp. timelike) Jordan Ivanov-Petrova for
every point $P$ of $M$. Similarly, we say that $(M,g)$ is {\it $k$-spacelike} (resp. {\it $k$-timelike})
Jordan Stanilov if $\mathcal{T}_P$ is $k$-spacelike (resp. $k$-timelike) Jordan Stanilov for
every point $P$ of $M$. In both contexts, note that the Jordan normal form is allowed to
vary with the point in question.

In the Riemannian setting $(p=0)$, the Jordan normal form is determined by the eigenvalue
structure and, as every $k$ plane is spacelike, we shall drop the qualifiers
`spacelike' and `Jordan'. This is not true in the higher signature context
which is why we focus on the Jordan normal form, i.e. the conjugacy class, instead of only
on the eigenvalue structure.

\subsection{Ivanov-Petrova tensors and manifolds} One has the following result, which is due to Gilkey, Leahy, and
Sadofsky
\cite{refGLS} and Gilkey \cite{refGi99} in the Riemannian setting
$(p=0)$, which was generalized by Zhang
\cite{Za00,Za02} to the Lorentzian ($p=1$) setting, and which was extended by Stavrov \cite{St03} to the higher signature setting:
\begin{theorem}\label{thm-1.1}
Let $\mathcal{V}$ be a model space with constant spacelike rank $r$.
\begin{enumerate}\item If $p=0$ and if $q\ne 3,4,7$, then $r=2$.
\item If $q\ge11$, if $1\le p\le(q-6)/4$, and if the set $\{q,q+1,...,q+p\}$ does not contain a
power of $2$, then $r=2$.
\end{enumerate}
\end{theorem}

This result is important as one has the following classification result
\cite{refGLS,GiZa02a}. Let $g_V$ be a metric on a finite dimensional real vector space $V$. If $\phi$ is a self-adjoint linear map
of $V$, then we define an algebraic curvature tensor $R_\phi$ on $V$ by setting:
$$R_\phi(x,y)z:=g_V(\phi y,z)\phi x-g_V(\phi x,z)\phi y\,.$$

\begin{theorem}\label{thm-1.2}
 Let $q\ge5$. The following assertions are equivalent:
\begin{enumerate}
\item The model space $\mathcal{V}$ is spacelike rank $2$ Jordan Ivanov-Petrova.
\item There exists $C\ne0$ and a self-adjoint map $\phi$ of $V$ so
$R=cR_\phi$ where one of the following 3 conditions on $\phi$ holds:
\begin{enumerate}
\item $\phi$ is an isometry of $(V,g)$, i.e. $g(\phi x,\phi y)=g(x,y)\ \forall\ x,y\in V$.
\item $\phi$ is a para-isometry of $(V,g)$, i.e. $g(\phi x,\phi y)=-g(x,y)\ \forall\ x,y\in
V$.
\item $\phi^2=0$ and $\ker\phi$ contains no spacelike vectors.
\end{enumerate}
\end{enumerate}
\end{theorem}
 
In the metric setting, one has \cite{refGLS,GiZa02b,IvPe}:

\begin{theorem}\label{thm-1.3} Let $(M,g)$ be a connected spacelike Jordan
Ivanov-Petrova pseudo-Riemannian manifold of signature $(p,q)$. Assume either $(p,q)=(0,4)$ or that $q\ge5$. Assume that
$R(\pi)$ is not nilpotent for at least one spacelike $2$ plane in $TM$ and that $\mathcal{T}_P$ has spacelike rank $2$ for
all $P\in M$. Let
$R=cR_\phi$ be as in Theorem
\ref{thm-1.2} where $\phi=\phi(P)$. Then $\phi$ is an isometry, $\phi^2=\id$, and one of the following $2$ cases
holds:
\begin{enumerate}
\item $\phi=\pm\id$ and $(M,g)$ has constant sectional curvature.
\item $(M,g)=(I\times N,\varepsilon dt^2+f(t)ds^2_N)$ where $(I,dt^2)$ is an open interval
in
$\mathbb{R}$, where
$(N,ds^2_N)$ has constant sectional curvature $K$, where $\varepsilon=\pm1$, and where the
warping function $f(t):=\varepsilon Kt^2+At+B$ for $A^2-4\varepsilon KB\ne0$.
\end{enumerate}
\end{theorem}

Theorems \ref{thm-1.1}, \ref{thm-1.2}, and \ref{thm-1.3} complete the classification of
Ivanov-Petrova manifolds in the Riemannian setting for $m\ne 3,4,7$; the work  of Ivanov and
Petrova \cite{IvPe} uses entirely different methods and shows that Theorem \ref{thm-1.3} holds if $m=4$. The case
$m=3$ is exceptional and the case $m=7$ is open.

In Section \ref{Sect-2}, we will use Theorem \ref{thm-1.3} to establish the following result:

\begin{theorem}\label{thm-1.4} Let $(M,g)$ be a connected spacelike Jordan
Ivanov-Petrova pseudo-Riemannian manifold of signature $(p,q)$. Assume either that $(p,q)=(0,4)$ or that $q\ge5$. Assume
that
$R(\pi)$ is not nilpotent for at least one spacelike $2$ plane in $TM$ and that $\mathcal{T}_P$ has spacelike rank $2$ for
all $P\in M$. Then
\begin{enumerate}
\item $(M,g)$ is $k$-spacelike Jordan Stanilov for any $2\le k\le q$.
\item $(M,g)$ is $k$-timelike Jordan Stanilov for any $2\le k\le p$.
\end{enumerate}
\end{theorem}

We will also establish the following partial converse in the Riemannian setting.

\begin{theorem}\label{thm-1.5} Let $(M,g)$ be a connected Riemannian manifold of dimension $m$ where $m\ne 3,7$.
If $(M,g)$ is $2$-Stanilov, then $(M,g)$ is Ivanov-Petrova with constant spacelike rank $2$.
\end{theorem}

Theorem \ref{thm-1.5} is false in the higher signature context. In Section \ref{Sect-3}
we will discuss a family of manifolds that arise as hypersurfaces in flat space. Let $p\ge2$ and let $(\vec x,\vec y)$ be
coordinates on $\mathbb{R}^{2p}$ where $\vec x=(x_1,...,x_p)$ and $\vec y=(y_1,...,y_p)$. Let $f=f(\vec x)$ be a smooth function
on a connected open subset $\mathcal{O}\subset\mathbb{R}^p$. Let $H=(H_{ij})$ be the Hessian where
$H_{ij}:=\partial_i^x\partial_j^xf$. We define a metric
$g_f$ on $M_f:=\mathcal{O}\times\mathbb{R}^p$ of neutral signature $(p,p)$ by setting:
\begin{equation}\label{eqn-1.a}
g_f(\partial_i^x,\partial_j^x):=\partial_i^xf\cdot\partial_j^xf,\quad
g_f(\partial_i^x,\partial_j^y)=g_f(\partial_j^y,\partial_i^x)=\delta_{ij},\quad
g_f(\partial_i^y,\partial_j^y)=0\,.
\end{equation}
We will establish the following result in Section \ref{Sect-3}; these manifolds were
first introduced in
\cite{GIZ03}.
\begin{theorem}\label{thm-1.6} Assume $\Rank(H)\ge2$. Then $(M_f,g_f)$ is:
\begin{enumerate}
\item spacelike Jordan Ivanov-Petrova if and only if $\det(H)$ is never zero.
\item timelike Jordan Ivanov-Petrova if and only if $\det(H)$ is never zero.
\item $k$-spacelike and $k$-timelike Jordan Stanilov for $2\le k\le p$ for any $f$.
\end{enumerate}
\end{theorem}

Let $s\ge2$. In Section \ref{Sect-4}, we exhibit family of manifolds of signature $(2s,s)$ which are spacelike Jordan
Ivanov-Petrova but not timelike Jordan Ivanov-Petrova. Thus the notions spacelike and timelike are distinct. This family provides
the first example of spacelike Ivanov-Petrova manifolds of spacelike rank $4$. The manifolds will be $k$-spacelike Jordan Stanilov
for all admissible $k$; they will be $k$-timelike Jordan Stanilov only for $k=2s$.

Let $(\vec u,\vec t,\vec v)$ be coordinates on $\mathbb{R}^{3s}$ where $\vec u=(u_1,...,u_s)$, $\vec t=(t_1,...,t_s)$, and $\vec
v=(v_1,...,v_s)$. We define a metric of signature $(2s,s)$ on $M_{3s}:=\mathbb{R}^{3s}$ by setting:
\begin{equation}\label{eqn-1.b}
\begin{array}{ll}
g_{3s}(\partial_i^u,\partial_j^u)=-2\delta_{ij}\sum_{1\le k\le s} u_kt_k,&
g_{3s}(\partial_i^u,\partial_j^v)=g_{3s}(\partial_j^v,\partial_i^u)=\delta_{ij},\\
g_{3s}(\partial_i^u,\partial_j^t)=g_{3s}(\partial_j^t,\partial_i^u)=0,&
g_{3s}(\partial_i^t,\partial_j^t)=-\delta_{ij},\vphantom{\vrule height 11pt}\\
g_{3s}(\partial_i^t,\partial_j^v)=g_{3s}(\partial_j^v,\partial_i^t)=0,&
g_{3s}(\partial_i^v,\partial_j^v)=0\,.\vphantom{\vrule height 11pt}
\end{array}
\end{equation}

One says that $(M,g)$ is {\it curvature homogeneous} if there exists a model $\mathcal{V}$ and isomorphisms
$\psi_P:\mathcal{T}_P\rightarrow\mathcal{V}$ for all $P\in M$; see \cite{KTV91,KTV92} for further details.
These manifolds were first introduced in \cite{GiNi03,GiNi03b} to provide examples of curvature homogeneous spacelike Osserman
manifolds which are not locally homogeneous and where the Jacobi operator was nilpotent of order $3$.

\begin{theorem}\label{thm-1.7}
The manifolds $(M_{3s},g_{3s})$ are
\begin{enumerate}
\item Spacelike rank $4$ Jordan Ivanov-Petrova.
\item Not timelike Jordan Ivanov-Petrova.
\item $k$-Spacelike Jordan Stanilov for $2\le k\le s$.
\item $k$-Timelike Jordan Stanilov if and only if $k=2s$.
\end{enumerate}
\end{theorem}

Here is a brief guide to this paper. In Section \ref{Sect-2}, we establish Theorems \ref{thm-1.4} and \ref{thm-1.5}. In
Section \ref{Sect-3}, we review results of \cite{GIZ03} to sketch the proof of Theorem \ref{thm-1.6}. Section \ref{Sect-4}
comprises the body of this paper. We first determine the curvature tensor of the metric in question. Then we show the space is
curvature homogeneous and determine the model space. We complete the proof of Theorem \ref{thm-1.7} by establishing the
corresponding assertions for the model space.

\section{Relationships between the Stanilov and Ivanov-Petrova condition}\label{Sect-2}

\begin{proof}[Proof of Theorem \ref{thm-1.4}] Let $(M,g)$ be a connected spacelike Jordan Ivanov-Petrova pseudo-Riemannian
manifold of signature $(p,q)$ where $q\ge5$ or where $(p,q)=(0,4)$. Assume that $R(\pi)$ is not nilpotent for at least one
spacelike
$2$ plane in
$TM$ and that $\mathcal{T}_P$ has spacelike rank $2$ for all $P\in M$. We may then use Theorem \ref{thm-1.3} to see $R=cR_\phi$
where
$\phi^2=\id$ and $\phi$ is an isometry. If $\{e_1,e_2\}$ is an orthonormal basis for an oriented spacelike $2$ plane
$\pi$, then
$R(e_1,e_2)x=c\{g(\phi e_2,x)\phi e_1-g(\phi e_1,x)\phi e_2\}$. Let $x\perp\pi$. Then
\begin{equation}\label{eqn-2.x}
\begin{array}{lll}
  R(\pi):\phi e_2\rightarrow c\phi e_1,&
  R(\pi):\phi e_1\rightarrow-c\phi e_2,&
  R(\pi):x\rightarrow 0,\\
  R(\pi)^2:\phi e_2\rightarrow -c^2\phi e_2,&
  R(\pi)^2:\phi e_1\rightarrow -c^2\phi e_1,&
  R(\pi)^2:x\rightarrow 0\,.\vphantom{\vrule height 11pt}
\end{array}\end{equation}
Let $\rho_\pi$ be orthogonal projection on $\pi$. By Display (\ref{eqn-2.x}),
\begin{equation}\label{eqn-2.y}
R(\pi)^2=-c^2\rho_{\phi\pi}\,.
\end{equation}
If $\pi_k\in\Gr_{k,+}(M,g)$, then we can use Equation (\ref{eqn-2.y}) to obtain
$$\Theta(\pi_k)=-(k-1)c^2\rho_{\phi\pi_k}\,.$$
This shows $\Theta(\pi_k)$ has constant Jordan normal form and hence $(M,g)$ is $k$-spacelike Jordan Stanilov. The argument that
$(M,g)$ is $k$-timelike Jordan Stanilov is essentially the same modulo an appropriate change of signs and thus is omitted. 
\end{proof}

\begin{proof}[Proof of Theorem \ref{thm-1.5}]
Let $(M,g)$ be a Riemannian manifold which is $2$-Stanilov. Let $\{\lambda_i(\pi)\}$ be the eigenvalues of $\mathcal{R}(\pi)$ for
$\pi\in\OGr_{2,+}(T_PM)$. Then $\{\lambda_i^2(\pi)\}$ are the eigenvalues of $\Theta(\pi)$. Since these eigenvalues are
independent of $\pi$, since $\OGr_{2,+}(T_PM)$ is connected, and since the eigenvalues vary continously, we may conclude
that $\mathcal{R}(\pi)$ also has constant eigenvalues. Since the Jordan normal form is determined by the eigenvalue structure in
the postive definite setting, we can conclude that $(M,g)$ is Ivanov-Petrova.
\end{proof}

\section{Stanilov manifolds of neutral signature $(p,p)$}\label{Sect-3} 
The following manifolds were first
introduced in \cite{GIZ03} and we follow the discussion there to see the metric $g_f$ of Equation (\ref{eqn-1.a}) is a
hypersurface metric. Let
$\{U_1,...,U_p,V_1,...,V_p,W\}$ be a basis for
$\mathbb{R}^{2p+1}$ where
$p\ge2$. Introduce a non-degenerate inner product $\langle\cdot,\cdot\rangle$ on $\mathbb{R}^{2p+1}$ by defining:
\begin{eqnarray*}
&&\langle U_i,V_j\rangle=\langle V_j,U_i\rangle=\delta_{ij},\quad\langle W,W\rangle=1,\\
&&\langle U_i,U_j\rangle=\langle V_i,V_j\rangle=\langle U_i,W\rangle=\langle W,U_i\rangle
=\langle V_i,W\rangle=\langle W,V_i\rangle=0\,.
\end{eqnarray*}

Introduce coordinates
$(\vec x,\vec y)$ on
$\mathbb{R}^{2p}$ where
$\vec x=(x_1,...,x_p)$ and $\vec y=(y_1,...,y_p)$. Let $f=f(\vec x)$ be a smooth real valued function on
$\mathcal{O}\subset\mathbb{R}^p$. Define an embedding of $M_f:=\mathcal{O}\times\mathbb{R}^p$ in $\mathbb{R}^{2p+1}$ by setting
$$\Psi_f(\vec x,\vec y)=\textstyle\sum_{1\le i\le s}(x_iU_i+y_iV_i)+f(\vec x)W\,.$$
Let $g_f$ be the induced hypersurface pseudo-Riemannian metric on $M_f$;
\begin{eqnarray*}
&&g_f(\partial_i^x,\partial_j^x)=\partial_i^xf\cdot\partial_j^xf,\quad
  g_f(\partial_i^x,\partial_j^y)=g_f(\partial_j^y,\partial_i^x)=\delta_{ij},\quad
  g_f(\partial_i^y,\partial_j^y)=0\,.
\end{eqnarray*}
Let $H_{ij}:=\partial_i^x\partial_j^xf$ be the Hessian, let $L$ be the second fundamental form, and let $S$ be the shape
operator:
\begin{equation}\label{eqn-3.a}
\begin{array}{lll}
L(\partial_i^x,\partial_j^x)=H_{ij},&
  L(\partial_i^x,\partial_j^y)=L(\partial_j^y,\partial_i^x)=0,&
  L(\partial_i^y,\partial_j^y)=0,\\
g_f(S(\cdot),\cdot)=L(\cdot,\cdot),&S(\partial_i^x)=\textstyle\sum_jH_{ij}\partial_j^y,&
  S(\partial_i^y)=0\,.\vphantom{\vrule height 11pt}
\end{array}
\end{equation}
We then have
\begin{equation}\label{eqn-3.b}
\begin{array}{l}
R(Z_1,Z_2,Z_3,Z_4)=L(Z_1,Z_4)L(Z_2,Z_3)-L(Z_1,Z_3)L(Z_2,Z_4),\\
R(Z_1,Z_2)Z_3=g_f(S(Z_2),Z_3)S(Z_1)-g_f(S(Z_1),Z_3)S(Z_2)\,.\vphantom{\vrule height 11pt}
\end{array}\end{equation}
Let $\ell:=\Rank(H)=\dim\Range\{S\}$; by assumption $2\le\ell\le p$.

\begin{proof}[Proof of theorem \ref{thm-1.6}] Let $\mathcal{X}:=\Span\{\partial_i^x\}$ and
$\mathcal{Y}:=\Span\{\partial_i^y\}$. If
$Z_1$ and
$Z_2$ are arbitrary tangent vectors, then we may use Displays (\ref{eqn-3.a}) and (\ref{eqn-3.b}) to see that
$$R(Z_1,Z_2)\mathcal{X}\subset\mathcal{Y}\quad\text{and}\quad R(Z_1,Z_2)\mathcal{Y}=0\,.$$
Consequently $R(Z_1,Z_2)^2=0$ and thus $\Theta(\pi)=0$ for any spacelike or timelike $k$ plane $\pi$. Thus, trivially, $(M_f,g_f)$
is $k$-spacelike and $k$-timelike Jordan Stanilov for any admissible $k$.

Since $\mathcal{R}(\pi)^2=0$ for any oriented spacelike (resp. timelike) $2$-plane, the Jordan normal form of
$\mathcal{R}(\pi)$ is determined by $\Rank\{\mathcal{R}(\pi)\}$. Let $\{Z_1,Z_2\}$ be an orthonormal basis for $\pi$. Expand
$Z_\mu=X_\mu+Y_\mu$ for $X_\mu\in\mathcal{X}$ and $Y_\mu\in\mathcal{Y}$. Then $R(\pi)=R(X_1,X_2)$ as well. Note that
$\Rank(S)=\Rank(S)|_{\mathcal{X}}$. Since
$\pi$ is spacelike (resp. timelike),
$X_1$ and $X_2$ are linearly independent. We have
$$\Rank\{\mathcal{R}(\pi)\}=\left\{\begin{array}{l}
0\text{ if }S(X_1)\text{ and }S(X_2)\text{ are linearly dependent},\\
2\text{ if }S(X_1)\text{ and }S(X_2)\text{ are linearly independent}\,.
\end{array}\right.$$
If $\Rank(S)=p$, then $\{S(X_1),S(X_2)\}$ is a linearly independent set and thus $\Rank\{\mathcal{R}(\pi)\}=2$; consequently
$(M_f,g_f)$ is spacelike (resp. timelike) Jordan Ivanov-Petrova.
On the other hand, if $2\le\Rank(S)\le p-1$, then we may choose spacelike (resp. timelike) $2$ planes $\pi_1$ and $\pi_2$ so
$\Rank\{\mathcal{R}(\pi_1)\}=2$ and $\Rank\{\mathcal{R}(\pi_2)\}=0$ and thus $(M_f,g_f)$ is not spacelike (resp. timelike) Jordan
Ivanov-Petrova.
\end{proof}

\section{$k$-Stanilov manifolds in signature $(2s,s)$}\label{Sect-4}

\subsection{The curvature tensor of the manifolds $(M_{3s},g_{3s})$}
We adopt the notation of Display (\ref{eqn-1.b}). We begin our study of the manifold $(M_{3s},g_{3s})$ by showing:

\begin{lemma}\label{lem-4.1} Let $R_{3s}$ be the curvature tensor of
the pseudo-Riemannian manifold $(M_{3s},g_{3s})$ defined in Display (\ref{eqn-1.b}). Then the non-zero entries in
$R_{3s}$ are, up to the usual $\mathbb{Z}_2$ symmetries of Equation (\ref{eqn-1.xx}), given by:
$$R_{3s}(\partial_i^u,\partial_j^u,\partial_j^u,\partial_i^u)
   =|u|^2\quad\text{and}\quad
R_{3s}(\partial_i^u,\partial_j^u,\partial_j^u,\partial_i^t)=1\,.$$
\end{lemma}

\begin{proof} Let $i\ne j$. The non-zero Christoffel symbols of the second kind are given by:
$$\begin{array}{ll}
g_{3s}(\nabla_{\partial_i^u}\partial_i^u,\partial_i^u)=-t_i,& \\
g_{3s}(\nabla_{\partial_i^u}\partial_i^u,\partial_i^t)=u_i,&
g_{3s}(\nabla_{\partial_i^u}\partial_i^t,\partial_i^u)=g_{3s}(\nabla_{\partial_i^t}\partial_i^u,\partial_i^u)=-u_i,\\
g_{3s}(\nabla_{\partial_i^u}\partial_i^u,\partial_j^u)=t_j,&
g_{3s}(\nabla_{\partial_i^u}\partial_j^u,\partial_i^u)=
g_{3s}(\nabla_{\partial_j^u}\partial_i^u,\partial_i^u)=-t_j,\vphantom{\vrule height 11pt}\\
g_{3s}(\nabla_{\partial_i^u}\partial_i^u,\partial_j^t)=u_j,&
g_{3s}(\nabla_{\partial_i^u}\partial_j^t,\partial_i^u)=
g_{3s}(\nabla_{\partial_j^t}\partial_i^u,\partial_i^u)=-u_j\,.\vphantom{\vrule height 11pt}
\end{array}$$
We may then raise indices to see the non-zero covariant derivatives are given by:
\begin{eqnarray*}
&&\nabla_{\partial_i^u}\partial_i^u=-t_i\partial_i^v
  +\textstyle\sum_{k\ne i}t_k\partial_k^v-\textstyle\sum_{1\le k\le s}u_k\partial_k^t,\\
&&\nabla_{\partial_i^u}\partial_j^u=-t_j\partial_i^v-t_i\partial_j^v,\\
&&\nabla_{\partial_i^u}\partial_i^t=\nabla_{\partial_i^t}\partial_i^u=-u_i\partial_i^v,\quad\text{and}\\
&&\nabla_{\partial_i^u}\partial_j^t=\nabla_{\partial_j^t}\partial_i^u=-u_j\partial_i^v\,.
\end{eqnarray*}

We have $\nabla\partial_i^v=0$. Thus if at least one $z_\mu\in\{\partial_i^v\}$,
$R_{3s}(z_1,z_2,z_3,z_4)=0$. Similarly, if at least two of the $z_\mu$ belong to $\{\partial_i^t\}$, then
$R_{3s}(z_1,z_2,z_3,z_4)=0$.  Furthermore 
$R_{3s}(\partial_i^u,\partial_j^u,\partial_k^u,\star)=0$ if the indices $\{i,j,k\}$ are distinct. Finally,
$$
\nabla_{\partial_i^u}\nabla_{\partial_j^u}\partial_j^u=-\partial_i^t+|u|^2\partial_i^v\quad\text{and}\quad
\nabla_{\partial_j^u}\nabla_{\partial_i^u}\partial_j^u=0\,.
$$
The Lemma now follows.\end{proof}

\begin{definition}\label{defn-4.2}
\rm  Let $\{U_1,...,U_s,T_1,...,T_s,V_1,...,V_s\}$ be a basis for $\mathbb{R}^{3s}$ where
$s\ge2$. Let $\mathcal{V}_{3s}:=(\mathbb{R}^{3s},g_{3s},R_{3s})$
where the non-zero entries of the metric $g_{3s}$ and of the algebraic curvature tensor
$R_{3s}$, up to the usual $\mathbb{Z}_2$ symmetries, are
\begin{equation}\label{eqn-4.a}
\begin{array}{l}
g_{3s}(U_i,V_i)=g_{3s}(V_i,U_i)=1,\quad
g_{3s}(T_i,T_i)=-1,\quad\text{and}\\
R_{3s}(U_i,U_j,U_j,T_i)=1\quad\text{for}\quad i\ne j\,.\vphantom{\vrule height 11pt}
\end{array}\end{equation}
Set $Z_i^\pm:=U_i\pm\ffrac12V_i$. Then $\Span\{Z_i^+\}$ is a maximal
spacelike subspace of $\mathbb{R}^{3s}$ and $\Span\{T_i,Z_i^-\}$ is the complementary maximal timelike subspace. Thus
$\mathbb{R}^{3s}$ has signature $(2s,s)$. A basis $\BB=\{\tilde U_1,...,\tilde U_s,\tilde T_1,...,\tilde T_s,\tilde
V_1,...,\tilde V_s\}$ for
$\mathbb{R}^{3s}$ is said to be {\it normalized} if the relations given above in Display (\ref{eqn-4.a}) hold for $\BB$. 
\end{definition}

\begin{lemma}\label{lem-4.3} $(M_{3s},g_{3s})$ is curvature homogeneous with model space $\mathcal{V}_{3s}$.
\end{lemma}

\begin{proof}Fix $P\in M_{3s}$. Let constants $\varepsilon_i$ and $\varrho_i$ be given. We
define a new basis for
$T_PM$ by setting:
$$U_i:=\partial_i^u+\varepsilon_{i}\partial_i^t+\varrho_{i}\partial_i^v,\quad
 T_i:=\partial_i^t+\varepsilon_{i}\partial_i^v,\quad\text{and}\quad V_i:=\partial_i^v\,.$$
Let $i\ne j$. Since $g_{3s}(U_i,T_i)=\varepsilon_{i}-\varepsilon_{i}=0$, the possibly non-zero entries of $g_{3s}$ and $R_{3s}$
are, up to the usual
$\mathbb{Z}_2$ symmetries, given by
\begin{eqnarray*}
&&g_{3s}(U_i,U_i)=g_{3s}(\partial_i^u,\partial_i^u)-\varepsilon_{i}^2+2\varrho_{i},\\
&&g_{3s}(T_i,T_i)=-1,\quad g_{3s}(U_i,V_i)=1,\\
&&R_{3s}(U_i,U_j,U_j,T_i)=1,\quad\text{and}\\
&&R_{3s}(U_i,U_j,U_j,U_i)=|u|^2+2\varepsilon_i+2\varepsilon_j\,.
\end{eqnarray*}
We set
$$\varepsilon_{i}:=-\ffrac14|u|^2
 \quad\text{and}\quad\varrho_{i}:=\ffrac12\varepsilon_{i}^2-\ffrac12g_{3s}(\partial_i^u,\partial_i^u)\,.
$$
This ensures that $g_{3s}(U_i,U_i)=0$ and
$R_{3s}(U_i,U_j,U_j,U_i)=0$ and establishes the existence of a basis with the normalizations of Definition \ref{defn-4.2}.
\end{proof}

Lemma \ref{lem-4.3} shows that the manifold $(M_{3s},g_{3s})$ is curvature homogeneous; the work of \cite{GiNi03b}
shows it is not locally homogeneous. We shall prove Theorem \ref{thm-1.7} by establishing the corresponding assertions for the
model space $\mathcal{V}_{3s}$.

\subsection{The skew-symmetric curvature operator} Theorem \ref{thm-1.7} (1,2) will follow from the following result concerning
the model space $\mathcal{V}_{3s}$.

\begin{lemma}\label{lem-4.4}
 Let $\mathcal{R}_{3s}$ be the skew-symmetric curvature operator defined by $R_{3s}$.\begin{enumerate}
\item If $\pi$ is an oriented spacelike $2$ plane, then
$\Rank\{\mathcal{R}_{3s}(\pi)\}=4$,\newline
$\Rank\{\mathcal{R}_{3s}(\pi)^2\}=2$, and $\mathcal{R}_{3s}(\pi)^3=0$.
\item The model space $\mathcal{V}_{3s}$ is spacelike rank $4$ Jordan Ivanov-Petrova.
\item The model space $\mathcal{V}_{3s}$ is not timelike Jordan Ivanov-Petrova.
\end{enumerate}\end{lemma}

\begin{proof}There is an additional useful symmetry which plays a crucial role. Let 
$$O(s):=\{\xi=(\xi_{ij}):\textstyle\sum_{1\le i\le s}\xi_{ij}\xi_{ik}=\delta_{jk}\}\subset\mathbb{M}_{3s}(\mathbb{R})$$
 be the standard orthogonal group of $s\times s$ real matrices. Define a diagonal action of
$O(s)$ on $\mathbb{R}^{3s}$ which preserves the structures $g_{3s}$ and $R_{3s}$ by setting:
\begin{equation}\label{eqn-4.b}
\xi:U_i\rightarrow\textstyle\sum_j\xi_{ij}U_j,\quad
 \xi:T_i\rightarrow\textstyle\sum_j\xi_{ij}T_j,\quad\text{and}\quad
 \xi:V_i\rightarrow\textstyle\sum_j\xi_{ij}V_j\,.
\end{equation}

Let $\pi$ be an oriented spacelike $2$ plane. By applying a symmetry of the form given in Equation (\ref{eqn-4.b}), we may suppose
that
$\pi=\operatorname{Span}\{X_1,X_2\}$ where
\begin{eqnarray*}
&&X_1=U_1+\textstyle\sum_{1\le i\le s}\{b_{1i}T_i+c_{1i}V_i\}\quad\text{and}\\
&&X_2=U_2+\textstyle\sum_{1\le i\le s}\{b_{2i}T_i+c_{2i}V_i\}\,.
\end{eqnarray*}
Let $c:=g_{3s}(X_1,X_1)g_{3s}(X_2,X_2)-g_{3s}(X_1,X_2)^2>0$ and let $\Xi:=R_{3s}(X_1,X_2)$. We then have
$\mathcal{R}_{3s}(\pi)=\frac1{\sqrt c}\Xi$. Let $i\ge3$. There exist real numbers $\varepsilon_{ij}=\varepsilon_{ij}(b,c)$ and
$\varrho_{ij}=\rho_{ij}(b,c)$, which play no role in the subsequent development, so that
\begin{equation}\label{eqn-4.c}
\begin{array}{lll}
\Xi:U_1\rightarrow T_2+\textstyle\sum_{1\le k\le s}\varepsilon_{1k}V_k,&
\Xi:T_1\rightarrow-V_2,&
\Xi:V_1\rightarrow0,\\
\Xi:U_2\rightarrow-T_1+\textstyle\sum_{1\le k\le s}\varepsilon_{2k}V_k,&
\Xi:T_2\rightarrow V_1,&
\Xi:V_2\rightarrow 0,\\
\Xi:U_i\rightarrow\varrho_{i1}V_1+\varrho_{i2}V_2,&
\Xi:T_i\rightarrow0,&
\Xi:V_i\rightarrow0\,.
\end{array}\end{equation}
Assertion (1) now follows; Assertion (2) follows from Assertion (1). Let
$$\pi_1:=\Span\{T_1,T_2\}\quad\text{and}\quad\pi_2:=\Span\{Z_1^-,Z_2^-\}$$
be timelike $2$ planes with
$\mathcal{R}_{3s}(\pi_1)=0$ and
$\mathcal{R}_{3s}(\pi_2)\ne0$. Assertion (3) follows. \end{proof}

\subsection{The higher order curvature operator of $\mathcal{V}_{3s}$}
Define a positive semi-definite bilinear form $\tilde g$ on $\mathbb{R}^{3s}$ by setting
$$\begin{array}{lll}
\tilde g(U_i,U_j)=\delta_{ij},&\tilde g(U_i,V_j)=\tilde g(V_j,U_i)=0,&\tilde g(U_i,T_j)=\tilde g(T_j,U_i)=0,\\
\tilde g(T_i,T_j)=0,&\tilde g(T_i,V_j)=\tilde g(V_j,T_i)=0,&\tilde g(V_i,V_j)=0\,.\vphantom{\vrule height 11pt}
\end{array}$$
This inner product is invariant under the action of $O(s)$ described in Equation (\ref{eqn-4.b}).
If $\pi$ is a linear subspace of $\mathbb{R}^{3s}$, set
$$\ell(\pi):=\Rank\{\tilde g|_\pi\}\,.
$$
 We complete the proof of Theorem \ref{thm-1.7} by showing:
\begin{lemma}\label{lem-4.5}
 Let $\Theta_{3s}$ be the higher order curvature operator defined by $\mathcal{V}_{3s}$.\begin{enumerate}
\item If $\pi$ is a spacelike $k$-plane, then
$\Rank\{\Theta_{3s}(\pi)\}=k$ and $\Theta_{3s}(\pi)^2=0$.
\item If $2\le k\le s$, then $\mathcal{V}_{3s}$ is $k$-spacelike Jordan Stanilov.
\item Let $\pi$ be a timelike $2s$ plane. If $\ell(\pi)\ge 2$, then $\Rank\{\Theta_{3s}(\pi)\}=\ell(\pi)$ and
$\Theta_{3s}(\pi)^2=0$. If $\ell(\pi)\le 1$, then $\Theta_{3s}(\pi)=0$.
\item  $\mathcal{V}_{3s}$ is $k$-timelike Jordan Stanilov if and only if $k=2s$.
\end{enumerate}\end{lemma}

\begin{proof} Fix $2\le k\le s$. Let indices $\alpha$ and $\beta$ range from $1$ through $k$. Let $\pi$ be a spacelike $k$-plane
in
$\mathbb{R}^{3s}$. We diagonalize the quadratic form $\tilde g|_\pi$ with respect to the positive definite quadratic form
$g|_\pi$ to choose an orthonormal basis $\{X_\alpha\}$ for $\pi$ so $\tilde g(X_\alpha,X_\beta)=a_\alpha
a_\beta\delta_{\alpha,\beta}$ where $a_\alpha>0$. By replacing $\pi$ by $\xi\cdot\pi$ for an appropriately chosen symmetry $\xi$
in
$O(s)$, we may assume without loss of generality that
$$X_\alpha=a_\alpha U_\alpha+\textstyle\sum_{1\le i\le s}\{b_{\alpha i}T_i+c_{\alpha i}V_i\}\quad\text{for}\quad1\le\alpha\le
k,$$ 
where the real numbers $b_{\alpha i}$ and $c_{\alpha i}$ play no role in the subsequent discussion.
Let $\alpha\ne\beta$. Let $1\le
i\le s$. We use Equation (\ref{eqn-4.c}) to see that:
\begin{eqnarray*}
&&R_{3s}(X_\alpha,X_\beta)^2U_i=\left\{\begin{array}{ll}
a_\alpha^2a_\beta^2V_i&\text{if }i=\alpha,\beta,\\
0&\text{otherwise}\,.
\end{array}\right.\\
&&R_{3s}(X_\alpha,X_\beta)^2T_i=0,\quad\text{and}\quad
  R_{3s}(X_\alpha,X_\beta)^2V_i=0\,.
\end{eqnarray*}
Since $\Theta_{3s}(\pi)=\textstyle\sum_{\alpha,\beta}R_{3s}(X_\alpha,X_\beta)^2
   =\sum_{\alpha,\beta}a_\alpha^2a_\beta^2R_{3s}(U_\alpha,U_\beta)^2$,
\begin{eqnarray*}
&&\Theta_{3s}(\pi)U_i=\left\{\begin{array}{lll}
\textstyle\sum_{\beta\ne i, 1\le\beta\le k}a_i^2a_\beta^2V_i&\text{if}&i\le k,\\
0&\text{if}&k+1\le i\le s,\end{array}\right.\\
&&\Theta_{3s}(\pi)T_i=0,\quad\text{and}\quad\Theta_{3s}(\pi)V_i=0\,.
\end{eqnarray*}
Assertion (1) now follows; Assertion (2) follows from Assertion (1).

Suppose that $\pi$ is a timelike $2s$ plane. We apply exactly the same diagonalization argument to see that, after replacing $\pi$
by
$\xi\cdot\pi$ for suitably chosen
$\xi\in O(s)$, we may assume without loss of generality there exists an orthonormal basis for
$\pi$ so
\begin{eqnarray*}
&&X_i=a_iU_i+\textstyle\sum_j\{b_{ij}T_j+c_{ij}V_j\}\quad\text{for}\quad1\le i\le \ell,\\
&&X_i=\textstyle\sum_j\{b_{ij}T_j+c_{ij}V_j\}\quad\text{for}\quad\ell+1\le i\le k\,,
\end{eqnarray*}
where $a_i>0$ for $1\le i\le\ell$. We then have
$$\Theta_{3s}(\pi)=\textstyle\sum_{1\le i,j\le\ell}a_i^2a_j^2R_{3s}(U_i,U_j)^2\,.$$

Since $\Theta_{3s}(\pi)^2=0$ for any $\pi$, the Jordan normal form is determined by $\Rank\{\Theta(\pi)\}$.
The argument given above shows 
$$\Rank\{\Theta_{3s}(\pi)\}=\left\{\begin{array}{lll}
\ell&\text{if}&\ell\ge2,\\
0&\text{if}&\ell<2\,.
\end{array}\right.$$
Assertion (3) follows. 

We use Assertion (3) to see that $\mathcal{V}_{3s}$ is $2s$ Jordan Stanilov. If $2\le k\le 2s-1$, then set:
\begin{eqnarray*}
\pi_1:&=&\left\{\begin{array}{ll}
\Span\{T_1,...,T_k\}&\text{if }2\le k\le s,\\
\Span\{T_1,...,T_s,Z_1^-,...,Z_{k-s}^-\}
\hphantom{.......}&\text{if }s<k<2s,\vphantom{\vrule height 12pt}
\end{array}\right.\\
\pi_2:&=&\left\{\begin{array}{ll}
\Span\{T_1,...,T_{k-2},Z_1^-,Z_2^-\}&\text{if }2\le k\le s,\\
\Span\{T_1,...,T_{s-1},Z_1^-,...,Z_{k+1-s}^-\}&\text{if }s<k<2s\,.\vphantom{\vrule height 12pt}
\end{array}\right.
\end{eqnarray*}
Then
\begin{eqnarray*}
&&\Rank\Theta_{3s}(\pi_1)=\left\{\begin{array}{ll}
0&\text{if }2\le k\le s+1,\\
k-s\hphantom{......}&\text{if }s+2\le k<2s,
\end{array}\right.\\
&\ne&\Rank\Theta_{3s}(\pi_2)=\left\{\begin{array}{ll}
2&\text{if }2\le k\le s+1,\\
k+1-s&\text{if }s+2\le k<2s\,.
\end{array}\right.
\end{eqnarray*}
This shows $\mathcal{V}_{3s}$ is not $k$-timelike Jordan Stanilov.
\end{proof}

\subsection{Remark:}\label{rmk-4.4} One can generalize the pseudo-Riemannian
manifold $(M_{3s},g_{3s})$ as follows. Let $\vec u:=(u_1,...,u_s)$, $\vec t:=(t_1,...,t_s)$, and $\vec v:=(v_1,...,v_s)$ give
coordinates $(\vec u,\vec t,\vec v)$ on $\mathbb{R}^{3s}$ for $s\ge2$. Let 
$F(\vec u):=f_1(u_1)+...+f_s(u_s)$
be a smooth function
on an open subset $\mathcal{O}\subset\mathbb{R}^s$.
Define a
pseudo-Riemannian metric
$g_F$ of signature
$(2s,s)$ on $M_F:=\mathcal{O}\times\mathbb{R}^{2s}$ whose non-zero components are given by:
$$
\begin{array}{l}
g_F(\partial_i^u,\partial_i^u)=-2F(\vec u)-2\textstyle\sum_iu_it_i,\\
 g_F(\partial_i^u,\partial_i^v)=g_F(\partial_i^v,\partial_i^u)=1,\\
 g_F(\partial_i^t,\partial_i^t)=-1\,.
\end{array}$$
It was shown in \cite{GiNi03b} that these spaces are curvature homogeneous with model space $\mathcal{V}_{3s}$ and thus
arguments given above show that the conclusions of Theorem
\ref{thm-1.7} apply to all of the manifolds in this family. These manifolds are not locally homogeneous for generic members of
the family.

\section*{Acknowledgments} Research of P. Gilkey partially supported by the
MPI (Leipzig). Research of S. Nik\v cevi\'c partially supported by the DAAD (Germany) and MM 1646 (Srbija).
Research of V. Videv partially supported by NSFI under contract MM809/98 (Bulgaria). The first two authors
wish to express their thanks to the Technical University of Berlin where much of the research reported here was conducted.
Finally, it is a pleasant task to thank Professor E. Garc\'{\i}a--R\'{\i}o for helpful discussions.


\begin{thebibliography}{AAA}

\bibitem{refDG} C. Dunn and P. B. Gilkey, 
{\it Curvature homogeneous pseudo-Riemannian manifolds which are not locally homogeneous}; 
math.DG/0306072.



\bibitem{refGi99} P. Gilkey,
 {\it Riemannian manifolds whose skew-symmetric curvature operator has constant eigenvalues II}, 
 {\bf Differential geometry and applications} (Brno, 1998), 73--87, Masaryk Univ., Brno, 1999.
ISBN: 80-210-2097-0

\bibitem{refGi02} ---, %%%P. Gilkey
{\bf Geometric Properties of Natural Operators Defined by the Riemann Curvature Tensor},
World Scientific ISBN 981-02-4752-4 (2002).

\bibitem{refGLS} P. Gilkey, J. V. Leahy, and H. Sadofsky,
{\it Riemannian manifolds whose skew-symmetric curvature operator has constant eigenvalues}, 
Indiana Univ. Math. J. \bf 48 \rm(1999), 615--634.

\bibitem{GIZ03} P. Gilkey, R. Ivanova and T. Zhang,
{\it Szabo Osserman IP Pseudo-Riemannian manifolds}, 
Publ. Math. Debrecen {\bf 62} (2003), 387--401; math.DG/0205085

\bibitem{GiNi03} P. Gilkey and S. Nik\v cevi\'c, 
{\it Nilpotent Spacelike Jorden Osserman pseudo-Riemannian manifolds}; 
math.DG/0302044.

\bibitem{GiNi03b} ---,
{\it Curvature homogeneous spacelike Jordan Osserman pseudo-Riemannian manifolds},
math.DG/0310024.

\bibitem{GiZa02a} P. Gilkey and T. Zhang,
{\it Algebraic curvature tensors whose skew-symmetric
curvature operator has constant rank 2},
Periodica Mathematica  Hungarica {\bf 44} (2002), 7--26

\bibitem{GiZa02b} ---, 
{\it Algebraic curvature tensors for indefinite metrics whose
  skew-symmetric curvature operator has constant Jordan normal form},
 Special issue for S. S. Chern. 
 Houston J. Math. {\bf 28} (2002), 311--328.

\bibitem{IvPe} S. Ivanov and I. Petrova,
{\it Riemannian manifold in which the skew-symmetric curvature operator has pointwise constant eigenvalues}, 
Geom. Dedicata {\bf 70} (1998), 269--282.

\bibitem{refRIGS} R. Ivanova and G. Stanilov, 
{\it A skew-symmetric curvature operator in Riemannian geometry},
{\bf Proceedings of the 2nd Gauss Symposium}. Conference A: Mathematics and Theoretical Physics (Munich, 1993), 391--395, 
Sympos. Gaussiana, de Gruyter, Berlin, 1995.

\bibitem{KTV91} O. Kowalski, F. Tricerri, and L. Vanhecke, 
{\it New examples of non-homogeneous Riemannian manifolds whose curvature tensor is that of a Riemannian symmetric space}, 
C. R. Acad. Sci. Paris S\'er. I Math. {\bf 311} (1990), 355--360.

\bibitem{KTV92} ---,
{\it Curvature homogeneous Riemannian manifolds},
 J. Math. Pures Appl. {\bf 71} (1992), 471--501.

\bibitem{Stx} G. Stanilov,
{\it Curvature operators based on the skew-symmetric curvature operator and their place in the Differential Geometry},
preprint (2000).

\bibitem{Sty} ---,
{\it Higher order skew-symmetric and symmetric curvature operators},
preprint (2003).

\bibitem{St03} I. Stavrov, {\bf Spectral geometry of the Riemann curvature tensor}, Ph. D.
Thesis, University of Oregon (2003).

\bibitem{TzVi99} J.Tzankov and V.Videv, 
{\it A Riemannian pointwise Stanilov manifolds of type (n,k)},
Abstracts of 4th International Conference on Geometry and Applications, Varna, 1999, 62--63.

\bibitem{Za00} T. Zhang, 
 {\it Manifolds with indefinite metrics whose skew-symmetric curvature operator has constant eigenvalues}
  Ph. D. thesis, University of Oregon (2000).

\bibitem{Za02} ---,
 {\it Applications of algebraic topology in bounding the rank of the skew-symmetric curvature operator},
 Topology Appl. {\bf 124} (2002), 9--24.
\end{thebibliography}
\end{document}